\documentclass[final, 1p, times]{elsarticle}
\usepackage{amsmath, amssymb, amsthm, amscd}
\usepackage{amsfonts}
\usepackage{graphicx}
\usepackage{color}
\usepackage{mathrsfs}
\usepackage{titletoc}
\usepackage{titlesec}
\usepackage{framed}

\usepackage{geometry}
\usepackage{setspace}

\newtheorem{theorem}{Theorem}

\newtheorem{lemma}[theorem]{Lemma}

\newcommand{\ra}{\rightarrow}
\newcommand{\p}{\partial}
\newcommand{\f}{\frac}

\renewcommand{\f}{\frac}

\newcommand{\be}{\begin{eqnarray}}
\newcommand{\ee}{\end{eqnarray}}
\newcommand{\bea}{\begin{eqnarray}}
\newcommand{\eea}{\end{eqnarray}}
\newcommand{\bna}{\begin{eqnarray*}}
\newcommand{\ena}{\end{eqnarray*}}

\renewcommand{\le}{\left}
\newcommand{\ri}{\right}

\usepackage{bbding}

\begin{document}

\begin{frontmatter}
\title{Topological degree for Chern-Simons Higgs models on finite graphs \tnoteref{sw}}

\author[ljy]{Jiayu Li}
\ead{jiayuli@ustc.edu.cn}

\author[sll]{Linlin Sun}
 \ead{sunlinlin@gxnu.edu.cn}

\author[yyy]{Yunyan Yang}
 \ead{yunyanyang@ruc.edu.cn}

\address[ljy]{School of Mathematical Sciences,
University of Science and Technology of China, Hefei 230026, China}
\address[sll]{School of Mathematics and Statistics, Guangxi Normal University, Guilin 541004, China}
 \address[yyy]{School of Mathematics,
Renmin University of China, Beijing 100872, China}

\tnotetext[sw]{This research is partly supported by the National Natural Science Foundation of China, Grant No.11721101 and Grant No. 12031017.}

\begin{abstract}
Let $(V,E)$ be a finite connected graph. We are concerned about the Chern-Simons Higgs model
$$\Delta u=\lambda e^u(e^u-1)+f, \eqno{(0.1)}$$
where $\Delta$ is the graph Laplacian, $\lambda$ is a real number and $f$ is a function on $V$. When $\lambda>0$ and
$f=4\pi\sum_{i=1}^N\delta_{p_i}$, $N\in\mathbb{N}$, $p_1,\cdots,p_N\in V$, the equation (0.1) was investigated by Huang, Lin, Yau (Commun. Math. Phys. 377 (2020) 613-621)
and Hou, Sun (Calc. Var. 61 (2022) 139) via the upper and lower solutions principle.
We now consider an arbitrary real number $\lambda$ and a general function $f$, whose integral mean is denoted by $\overline{f}$, and prove that when $\lambda\overline{f}<0$,  the equation $(0.1)$ has a solution; when $\lambda\overline{f}>0$,
there exist two critical numbers $\Lambda^\ast>0$ and $\Lambda_\ast<0$ such that if $\lambda\in(\Lambda^\ast,+\infty)\cup(-\infty,\Lambda_\ast)$,
then $(0.1)$ has at least two solutions, including one local minimum solution; if $\lambda\in(0,\Lambda^\ast)\cup(\Lambda_\ast,0)$, then $(0.1)$ has no solution; while if $\lambda=\Lambda^\ast$ or $\Lambda_\ast$, then $(0.1)$ has at least one solution.
Our method is calculating the topological degree and using the relation between the degree and the critical group of a related functional.
Similar method is also applied to the Chern-Simons Higgs system, and a partial result for the multiple solutions of the system is obtained.
\end{abstract}

\end{frontmatter}

\section{Introduction}
The Chern-Simons Higgs model, introduced by Hong, Kim, Pac \cite{HongKP} and Jackiw, Weinberg \cite{JackiwW},
has always attracted the attention of many mathematicians in the fields of geometry and physics, see for examples
\cite{Carffarali, Chae, DJLW1,DJLW2,LinPY,Nolasco,Nolasco-2,Tarantello,Wang}. Among many versions, the self-dual
Chern-Simons Higgs vortex equation on a flat 2-torus $\Sigma$ can be written as
\be\label{eeq-1}\Delta u=\f{4}{k^2}e^u(e^u-1)+4\pi\sum_{i=1}^{k_0}m_i\delta_{p_i},\ee
where $k>0$ is the Chern-Simons constant, $m_i\in\mathbb{N}$, $p_i\in\Sigma$, $i=1,\cdots,k_0$. The solution of the above equation is called
a vertex solution, each $p_i$ is called a vertex point, and $m_i$ stands for the multiplicity of $p_i$. From the view of physics, the vortex points
are closely related to the local maximum point of the magnetic flux in the Chern-Simons Higgs model. Let $u_0$ be a solution of
$$\le\{\begin{array}{lll}\Delta u_0=-\f{4\pi N}{|\Sigma|}+4\pi\sum_{i=1}^{k_0}m_i\delta_{p_i}\\[1.2ex]
\int_\Sigma u_0dv_g=0,\end{array}\ri.$$
where $N=\sum_{i=1}^{k_0}m_i$. Set $v=u-u_0$. Then (\ref{eeq-1}) can be written in a more favourable form
\be\label{eeq-2}\Delta v=\lambda he^v(he^v-1)+\f{4\pi N}{|\Sigma|},\ee
where $\lambda=\f{4}{k^2}$ and $h=e^{u_0}$ is a positive function on $\Sigma$.
A solution $v$ of (\ref{eeq-2}) is called of finite energy if $v\in W^{1,2}(\Sigma)$, a usual Sobolev space.
Indeed, it is known that the corresponding physical energy of the solution $v$ is finite if $u\in W^{1,2}(\Sigma)$.
Thus, solutions of finite energy are physically meaningful in (\ref{eeq-2}) and there have been
many existence results for $W^{1,2}(\Sigma)$ solutions of (\ref{eeq-2}), see \cite{Carffarali, DJLW-0,DJLW1,Nolasco,Nolasco-2,Tarantello,Wangmeng-2,Wangmeng} and the references therein. By using
the principle of upper and lower solutions, Caffarelli and Yang constructed a maximal solution.
In addition to the above references, \cite{DJLW2,lanli} also indicated
that the equation (\ref{eeq-2}) admits a variational structure.

Different from the theoretical significance on Riemann surfaces, the analysis on graphs is very important for applications, such as image processing, data mining,
network and so on. Among lots of directions, partial differential equations arising in geometry or physics are worth studying on graphs.
Various equations, including the heat equation \cite{Horn, Huang, Lin1, Lin2}, the Fokker-Planck and Schr\"odinger equations \cite{Chow1, Chow2}, have been studied by many mathematicians. In particular, Grigor'yan, Lin and Yang \cite{Gri1,Gri2,Gri3} studied the existence  of solutions for a series of nonlinear elliptic equations on graphs by using the variational methods. In this direction,
 Zhang, Zhao, Han and Shao \cite{Han2,Han,ZhangZhao} obtained  nontrivial solutions to nonlinear Schr\"odinger equations
 with potential wells. Similar problems on infinite metric graphs were studied by Akduman-Pankov \cite{Akduman-Pankov}.
 The Kazdan-Warner equation was extended by Keller-Schwarz \cite{Keller-Schwarz} to canonically compactifiable graphs.
 Semi-linear heat equations on locally finite graphs were studied by Ge, Jiang, Lin and Wu \cite{Ge-Jiang,Lin1,Lin2}.
 For other related works, we refer the readers to \cite{Ge-Proc,SunYH,Hua,Hua-2,Lin-Yang1,Lin-Yang2,LiuYang,LiuY,ref16,Zhang-Lin,Zhang-Lin2,zhuxb} and the references
 therein.

  To describe the Chern-Simons Higgs model in the graph setting, we introduce some notations.
  Let $(V,E)$ be a connected finite graph, where $V$ is the set of vertices and $E$ is the set of edges. Let
   $\mu :V\rightarrow (0,+\infty)$ and $\{w_{xy}:xy\in E\}$ be its measure and weights respectively. The weight $w_{xy}$ is always
   assumed to be positive and symmetric.
   The Laplacian of a function $u: V\rightarrow \mathbb{R}$ reads as
$$\Delta u(x)=\frac{1}{\mu (x)} \sum_{y\sim x} w_{xy}(u(y)-u(x)),$$
where $y\sim x$ means $y$ is adjacent to $x$, i.e. $xy\in E$. The gradient of $u$ is defined as
$$\nabla u(x)=\le(\sqrt{\f{w_{xy_1}}{2\mu(x)}}(u(y_1)-u(x)),\cdots,\sqrt{\f{w_{xy_{\ell_x}}}{2\mu(x)}}(u(y_{\ell_x})-u(x))\ri),$$
where $\{y_1,\cdots,y_{\ell_x}\}$ are all distinct points adjacent to $x$. Clearly, such an $\ell_x$ is unique and
$\nabla u(x)\in\mathbb{R}^{\ell_x}$. The integral of $u$ is given by
$$\int_{V}ud\mu=\sum_{x\in V}\mu(x)u(x).$$

Now we consider an analog of (\ref{eeq-2}) on a connected finite graph, namely
\be\label{eq-1}\Delta u=\lambda e^u(e^u-1)+f\quad{\rm in}\quad V,\ee
where $\lambda\in \mathbb{R}$, $f:V\ra\mathbb{R}$ is a function. It was proved by Huang, Lin and Yau \cite{LinYau} that
if $\lambda>0$ and $f=4\pi\sum_{i=1}^N\delta_{p_i}$, there exists a critical number $\lambda^\ast>0$ such that
(\ref{eq-1}) has a solution when $\lambda>\lambda^\ast$, while (\ref{eq-1}) has no solution when $0<\lambda<\lambda^\ast$.
The critical case $\lambda=\lambda^\ast$ was solved by Hou and Sun \cite{HouSun}, who proved that (\ref{eq-1}) has also a solution.
Such results are essentially based on the method of upper and lower solutions principle. This together with variational method may lead to
 existence results for other forms of Chern-Simons Higgs models, see Chao and Hou \cite{Hou-2}. Recently, a more delicate analysis was employed by Huang, Wang and Yang \cite{HuangWY} to get existence of solutions of the Chern-Simons Higgs system.

Topological degree theory is a powerful tool in studying partial differential equations in the Euclidean space
or Riemann surfaces, see for example Li
\cite{Liyanyan}. It was first used by
Sun and Wang \cite{SunW} to solve the Kazdan-Warner equation on finite graphs. Very recently, it was also
employed by Liu \cite{LiuY-3} to deal with the mean field equation. Our aim is to use this powerful tool to study the
Chern-Simons Higgs model. The first and most important step is to get a priori estimate for solutions, say

\begin{theorem}\label{prior-1}
Let $(V,E)$ be a connected finite graph with symmetric weights, i.e. $w_{xy}=w_{yx}$ for all $xy\in E$. Let $\sigma\in[0,1]$, $\lambda$ and $f$ satisfy
\be\label{cond-0}\Lambda^{-1}\leq |\lambda|\leq \Lambda,\,\,\,
\Lambda^{-1}\leq \le|\int_Vfd\mu\ri|\leq \Lambda,\,\,\,\|f\|_{L^\infty(V)}\leq \Lambda\ee
for some real number $\Lambda>0$. If $u$ is a solution of
\be\label{eq-hom}\Delta u=\lambda e^{u}(e^u-\sigma)+f\quad{\rm in}\quad V,\ee
then there exists a constant $C$, depending only on $\Lambda$ and the graph $V$, such that
$|u(x)|\leq C$ for all $x\in V$.
\end{theorem}

When $\sigma=1$, the equation (\ref{eq-hom}) is exactly (\ref{eq-1}). In the case $\lambda>0$ and $f=4\pi\sum_{i=1}^N\delta_{p_i}$,
where $p_1,\cdots,p_N\in V$ and $N\in\mathbb{N}$, let
 $\lambda^\ast$ be the critical number in \cite{LinYau}.
Then for any $\lambda_k>\lambda^\ast$ with $\lambda_k\ra \lambda^\ast$ as $k\ra\infty$, there exists a solution $u_{\lambda_k}$ of (\ref{eq-1})
with $\lambda=\lambda_k$, $k=1,2,\cdots$. It follows from Theorem \ref{prior-1}
that $(u_{\lambda_k})$ is uniformly bounded in $V$. Hence up to a subsequence,
$(u_{\lambda_k})$ uniformly converges to some $u^\ast$, which is a solution of (\ref{eq-1}) with $\lambda=\lambda^\ast$. This gives another
proof of a result of Hou and Sun \cite{HouSun}.

Denote $X=L^\infty(V)$ and define a map $F:X\ra X$ by
\be\label{map}F(u)=-\Delta u+\lambda e^{u}(e^u-1)+f.\ee
The second step is to calculate the topological degree of $F$ by using its homotopic invariance property.

\begin{theorem}\label{degree}
Let $(V,E)$ be a connected finite graph with symmetric weights, and $F:X\ra X$ be a map defined by (\ref{map}).
Suppose that $\lambda\int_V{f}d\mu\not=0$. Then
there exists a large number $R_0>0$ such that for all $R\geq R_0$,
$$\deg(F,B_R,0)=\le\{
\begin{array}{lll}
1&{\rm if}& \lambda>0,\,\int_V{f}d\mu<0\\[1.2ex]
0&{\rm if}& \lambda\int_V{f}d\mu>0\\[1.2ex]
-1&{\rm if}& \lambda<0,\,\int_V{f}d\mu>0,
\end{array}\ri.$$
where $B_R=\{u\in X:\|u\|_{L^\infty(V)}<R\}$ is a ball in $X$.
\end{theorem}

As an application of the above topological degree, our existence results for the Chern-Simons Higgs model read as follows:

\begin{theorem}\label{existence}
Let $(V,E)$ be a connected finite graph with symmetric weights. Then we have the following:\\
$(\mathsf{a})$ If $\lambda\int_Vfd\mu<0$, then the equation (\ref{eq-1}) has a solution;\\
 $(\mathsf{b})$ If $\lambda\int_Vfd\mu>0$, then two subcases are distinguished:
  $(i)$ $\int_Vfd\mu>0$. There exists a
real number $\Lambda^\ast>0$ such that when $\lambda>\Lambda^\ast$, (\ref{eq-1}) has at least two different solutions; when $0<\lambda<\Lambda^\ast$,
(\ref{eq-1}) has no solution; when $\lambda=\Lambda^\ast$, (\ref{eq-1}) has at least one solution;
$(ii)$ $\int_Vfd\mu<0$. There exists a real number $\Lambda_\ast<0$ such that when $\lambda<\Lambda_\ast$, (\ref{eq-1}) has at least two different solutions;
when $\Lambda_\ast<\lambda<0$,
(\ref{eq-1}) has no solution; when $\lambda=\Lambda_\ast$, (\ref{eq-1}) has at least one solution.
\end{theorem}

We remark that Case $(\mathsf{b})$ $(i)$ includes
$\lambda>0$ and $f=4\pi\sum_{i=1}^N\delta_{p_i}$ as a special case, which was studied in \cite{Hou-2,HouSun,LinYau,HuangWY}.
In the subcase $\lambda>\Lambda^\ast>0$ or $\lambda<\Lambda_\ast<0$, we shall construct a local minimum solution,
and then use the topological degree to obtain the existence of another solution.
Our arguments are essentially different
from those in \cite{Hou-2,HuangWY,LiuYang}. Note that a solution of (\ref{eq-1}) is a critical point of the functional
$J_\lambda:X\ra\mathbb{R}$ defined by
\be\label{funct-1}J_\lambda(u)=\f{1}{2}\int_V|\nabla u|^2d\mu+\f{\lambda}{2}\int_V(e^u-1)^2d\mu+\int_Vfud\mu.\ee
Here a local minimum solution of (\ref{eq-1}) means a local minimum critical point of $J_\lambda$.\\

Also we consider the Chern-Simons Higgs system
\be\label{system}
\le\{
\begin{array}{lll}
\Delta u=\lambda e^v(e^u-1)+f\\[1.2ex]
\Delta v=\lambda e^u(e^v-1)+g,
\end{array}\ri.
\ee
where $\lambda$ is a real number, and $f,g$ are functions on $V$. Similar to the single equation, we need also a priori estimate.

\begin{theorem}\label{system-apriori}
Let $(V,E)$ be a connected finite graph with symmetric weights.
Suppose that $\sigma\in[0,1]$, $\lambda,\eta$ are two positive real numbers, $f,g$ are two functions verifying that $\int_V{f}d\mu>0$
and $\int_V{g}d\mu>0$. If $(u,v)$ is a solution of the system
\be\label{system-8}
\le\{
\begin{array}{lll}
\Delta u=\lambda e^v(e^u-\sigma)+f\\[1.2ex]
\Delta v=\eta e^u(e^v-\sigma)+g,
\end{array}\ri.
\ee
 then there exists a constant $C$, depending only on
$\lambda,\eta,f,g$ and the graph $V$, such that
$$\|u\|_{L^\infty(V)}+\|v\|_{L^\infty(V)}\leq C.$$
\end{theorem}

To compute the topological degree, we define a map $\mathcal{F}:X\times X\ra X\times X$ by
\be\label{map-syst}\mathcal{F}(u,v)=(-\Delta u+\lambda e^{v}(e^u-1)+f,-\Delta v+\eta e^u(e^v-1)+g).\ee

\begin{theorem}\label{degree-system}
 Let $(V,E)$ be a connected finite graph with symmetric weights, and $\mathcal{F}$ be a map defined by (\ref{map-syst}). If $\lambda>0$, $\eta>0$, $\int_Vfd\mu>0$ and $\int_Vgd\mu>0$,
 then there exists a large number $R_0>0$ such that for all $R\geq R_0$,
 $$\deg(\mathcal{F},B_R,(0,0))=0,$$ where $B_R=\{(u,v)\in X\times X:\|u\|_{L^\infty(V)}+\|v\|_{L^\infty(V)}<R\}$ is a ball in
 $X\times X$.
\end{theorem}

Define a functional $\mathcal{J}_\lambda:X\times X\ra\mathbb{R}$ by
\be\label{functional}\mathcal{J}_\lambda(u,v)=\int_V\nabla u\nabla vd\mu+\lambda\int_V(e^u-1)(e^v-1)d\mu+\int_V(fv+gu)d\mu.\ee
Note that for all $(\phi,\psi)\in X\times X$,
\bea\nonumber
\langle\mathcal{J}_\lambda^\prime(u,v),(\phi,\psi)\rangle&=&\le.\f{d}{dt}\ri|_{t=0}\mathcal{J}(u+t\phi,v+t\psi)\\\label{derivative}
&=&\int_V\le\{\le(-\Delta v+\lambda e^u(e^v-1)+g\ri)\phi+\le(-\Delta u+\lambda e^v(e^u-1)+f\ri)\psi\ri\}d\mu.
\eea
Clearly $(u,v)$ is a critical point of $\mathcal{J}_\lambda$ if and only if it is a solution of the system (\ref{system}).
As a consequence of Theorem \ref{degree-system}, we have the following
\begin{theorem}\label{syst-thm}
Let $(V,E)$ be a connected finite graph with symmetric weights, $\lambda>0$, $\int_Vfd\mu>0$, $\int_Vgd\mu>0$, and
$\mathcal{J}_\lambda$ be a functional defined by (\ref{functional}).
If either $\mathcal{J}_\lambda$ has a non-degenerate critical point, or $\mathcal{J}_\lambda$ has a local minimum critical point,
then it must have another  critical point.
\end{theorem}

It should be remarked that Theorem \ref{syst-thm} gives another solution of (\ref{system}) under the condition that
$\mathcal{J}_\lambda$ has a non-degenerate or a local minimum critical point beforehand. So it is only a partial result
for the problem of multiple solutions of the system (\ref{system}).\\

The remaining part of this paper is organized as follows: In Section \ref{sec-apri}, we give a priori estimate for solutions of (\ref{eq-1})
(Theorem \ref{prior-1});  The topological degree of $F:X\ra X$ (Theorem \ref{degree}) was calculated in Section \ref{sec-Brouwer};
In Section \ref{sec-existence}, we prove the existence result (Theorem \ref{existence}); The priori estimate and existence of solutions of
the Chern-Simons Higgs system (Theorems \ref{system-apriori}-\ref{syst-thm}) are discussed in Section \ref{sec-system}.

\section{A priori estimate}\label{sec-apri}

In this section, we shall prove Theorem \ref{prior-1}. In order to provide readers with a clear understanding of the proof,
we demonstrate the entire process from simple cases to complex cases. Precisely the proof will be divided into several lemmas as below.\\

The first priori estimate is for fixed $\lambda$ and $f$.
\begin{lemma}\label{simp-1}
Suppose that $u$ is a solution of (\ref{eq-1}), where $\lambda\not= 0$ and $\int_Vfd\mu\not=0$. Then there exists a constant $C$,
depending only on $\lambda$, $f$ and the graph $V$, such that $|u(x)|\leq C$ for all $x\in V$.
\end{lemma}

\proof  If $u$ is a solution of (\ref{eq-1}), then integration by parts gives
\be\label{eq-2}0=\int_V\Delta ud\mu=\lambda\int_Ve^u(e^u-1)d\mu+\int_Vfd\mu.\ee
Firstly, we show that $u$ has a uniform upper bound. With no loss of generality, we may assume $\max_Vu>0$. For otherwise, $u$ has already
upper bound 0.
Observing
$$\le|\int_{u<0}e^u(e^u-1)d\mu\ri|\leq \f{1}{4}|V|,$$
we derive from (\ref{eq-2}) that
$$\int_{u\geq 0}e^u(e^u-1)d\mu\leq a:=\f{1}{4}|V|+\f{1}{|\lambda|}\le|\int_Vfd\mu\ri|.$$
This together with the fact
$$\int_{u\geq 0}e^u(e^u-1)d\mu=\sum_{x\in V,\,u(x)\geq 0}\mu(x)e^{u(x)}(e^{u(x)}-1)
\geq\mu_0e^{\max_Vu}(e^{\max_Vu}-1)$$
leads to
\be\label{upper}\max_V u\leq \log\f{1+\sqrt{1+4a/\mu_0}}{2},\ee
where $\mu_0=\min_{x\in V}\mu(x)>0$, since $V$ is finite.

Secondly, we prove that $u$ has also a uniform lower bound. To see this, in view of (\ref{eq-1}) and (\ref{upper}), we calculate
for any $x\in V$,
\bna
|\Delta u(x)|&\leq&|\lambda|\le|e^{u(x)}(e^{u(x)}-1)\ri|+|f(x)|\\
&\leq&|\lambda|(e^{2u(x)}+e^{u(x)})+|f(x)|\\
&\leq&|\lambda|\le(\f{(1+\sqrt{1+4a/\mu_0})^2}{4}+\f{1+\sqrt{1+4a/\mu_0}}{2}\ri)+\|f\|_{L^\infty(V)}\\
&=:&b.
\ena
Hence, there holds
\be\label{up-2-0}\|\Delta u\|_{L^\infty(V)}\leq b.\ee
We may assume $V=\{x_1,\cdots,x_\ell\}$, $u(x_{1})=\max_Vu$, $u(x_{\ell})=\min_Vu$, and without loss of generality
$x_1x_2, x_2x_3,\cdots, x_{\ell-1}x_\ell$ is the shortest path connecting $x_1$ and $x_\ell$. It follows that
\bea\nonumber
0\leq u(x_1)-u(x_\ell)&\leq&\sum_{j=1}^{\ell-1}|u(x_j)-u(x_{j+1})|\\\nonumber
&\leq&\f{\sqrt{\ell-1}}{\sqrt{w_0}}\le(\sum_{j=1}^{\ell-1}w_{x_jx_{j+1}}(u(x_j)-u(x_{j+1}))^2\ri)^{1/2}\\\label{ineq-0}
&\leq&\f{\sqrt{\ell-1}}{\sqrt{w_0}}\le(\int_V|\nabla u|^2d\mu\ri)^{1/2},
\eea
where $w_0=\min_{x\in V,\,y\sim x} w_{xy}>0$. Denoting $\overline{u}=\f{1}{|V|}\int_Vud\mu$, we obtain by integration by parts
\bna
\int_V|\nabla u|^2d\mu&=&-\int_V(u-\overline{u})\Delta ud\mu\\
&\leq&\le(\int_V(u-\overline{u})^2d\mu\ri)^{1/2}\le(\int_V(\Delta u)^2d\mu\ri)^{1/2}\\
&\leq&\le(\f{1}{\lambda_1}\int_V|\nabla u|^2d\mu\ri)^{1/2}\le(\int_V(\Delta u)^2d\mu\ri)^{1/2},
\ena
which gives
\be\label{ineq-1}\int_V|\nabla u|^2d\mu\leq \f{1}{\lambda_1}\int_V(\Delta u)^2d\mu\leq \f{1}{\lambda_1}
\|\Delta u\|_{L^\infty(V)}^2|V|,\ee
where $\lambda_1=\inf_{\overline{v}=0,\int_Vv^2d\mu=1}\int_V|\nabla v|^2d\mu>0$. Combining (\ref{ineq-0}) and (\ref{ineq-1}), we conclude
\be\label{equiv}\max_Vu-\min_V u\leq \sqrt{\f{(\ell-1)|V|}{w_0\lambda_1}}\|\Delta u\|_{L^\infty(V)}.\ee
We remark that (\ref{equiv}) holds for arbitrary function $u$, such an inequality was obtained by Sun and Wang \cite{SunW} by using
the equivalence of all norms in a finite dimensional vector space, and here we give an explicit constant instead of $C$.
The power of (\ref{equiv}) is evident. In view of (\ref{up-2-0}), we have
\be\label{scope}\max_Vu-\min_V u\leq c_0:=b\sqrt{\f{(\ell-1)|V|}{w_0\lambda_1}}.\ee
Coming back to (\ref{eq-2}), we have
\be\label{eq-3}\int_Ve^u(e^u-1)d\mu=c_1:=-\f{1}{\lambda}\int_Vfd\mu.\ee
By the assumptions $\lambda\not=0$ and $\int_Vfd\mu\not=0$, we know $c_1\not=0$.
Now we {\it claim} that
\be\label{A}\max_Vu> -A:=\log\min\le\{1,\f{|c_1|}{4|V|}\ri\}.\ee
For otherwise, $\max_Vu\leq -A$, which together with (\ref{eq-3}) implies
\bna
|c_1|&=&\le|\int_Ve^u(e^u-1)d\mu\ri|\\
&\leq&\int_V(e^{2u}+e^u)d\mu\\
&\leq&(e^{2\max_Vu}+e^{\max_Vu})|V|\\
&\leq& 2e^{-A}|V|\\
&<&\f{|c_1|}{2}.
\ena
This contradicts $c_1\not=0$, and thus confirms our claim (\ref{A}). Inserting (\ref{A}) into (\ref{scope}), we obtain
$$-A-c_0\leq \min_Vu\leq\max_V u\leq \log\f{1+\sqrt{1+4a/\mu_0}}{2},$$
as we desired. $\hfill\Box$\\

The second priori estimate is for the changing $\lambda$ and $f$.

\begin{lemma}\label{remark2}
Let $u$ be a solution of (\ref{eq-1}). If $\lambda$ and $f$ satisfy (\ref{cond-0}),
then there exists a constant $C$, depending only on $\Lambda$ and the graph $V$, such that
$|u(x)|\leq C$ for all $x\in V$.
\end{lemma}

\proof It suffices to modify the argument in the proof of Lemma \ref{simp-1}.

Similar to
(\ref{upper}), we first have the upper bound estimate
\be\label{ineq-2}\max_V u\leq \log\f{1+\sqrt{1+4a/\mu_0}}{2},\ee
where $\mu_0=\min_{x\in V}\mu(x)$ and $a=|V|+\Lambda^2$. Next, instead of (\ref{scope}), we have
\be\label{ineq-3}\max_Vu-\min_V u\leq c_0=b\sqrt{\f{(\ell-1)|V|}{w_0\lambda_1}},\ee
where $\lambda_1=\inf_{\overline{v}=0,\int_Vv^2d\mu=1}\int_V|\nabla v|^2d\mu$, $\ell$ denotes the number of all points of $V$,
$w_0=\min_{x\in V,\,y\sim x}w_{xy}$ and
$$b=\Lambda\le(\f{(1+\sqrt{1+4a/\mu_0})^2}{4}+\f{1+\sqrt{1+4a/\mu_0}}{2}+1\ri).$$
To proceed, we shall show
\be\label{ineq-4}\max_Vu> -A=\log\min\le\{1,\f{1}{4|V|\Lambda^2}\ri\}.\ee
Suppose not. We have $\max_Vu\leq -A$ and
\bna
\f{1}{\Lambda^2}\leq \le|\f{1}{\lambda}\int_Vfd\mu\ri|&=&\le|\int_Ve^u(e^u-1)d\mu\ri|\\
&\leq&\int_V(e^{2u}+e^u)d\mu\\
&\leq& 2e^{-A}|V|\\
&<&\f{1}{2\Lambda^2},
\ena
which is impossible. Thus (\ref{ineq-4}) holds. Combining (\ref{ineq-2}), (\ref{ineq-3}) and (\ref{ineq-4}), we get the desired result. $\hfill\Box$ \\

The third priori estimate is not only for changing $\lambda$ and $f$, but also for the changing parameter $\sigma$.

\begin{lemma}\label{apriori-2}
Let $\sigma\in [0,1]$, $\lambda$ and $f$ satisfy (\ref{cond-0}) for some real number $\Lambda>0$. If $u$ is a solution of
(\ref{eq-hom}), then there exists a constant $C$, depending only on $\Lambda$ and the graph $V$, such that $|u(x)|\leq C$ for
all $x\in V$.
\end{lemma}

\proof  If $u$ is a solution of (\ref{eq-hom}), then integration by parts gives
$$0=\int_V\Delta ud\mu=\lambda\int_Ve^u(e^u-\sigma)d\mu+\int_Vfd\mu.$$
Similar to (\ref{upper}), keeping in mind $\sigma\in[0,1]$, we first have the same upper bound estimate as (\ref{ineq-2}), namely
$$\max_V u\leq \log\f{1+\sqrt{1+4a/\mu_0}}{2},$$
where $\mu_0=\min_{x\in V}\mu(x)$ and $a=|V|+\Lambda^2$. Next, we have the same estimates as (\ref{ineq-3}) and (\ref{ineq-4}), which is
independent of the parameter $\sigma\in[0,1]$. In particular
$$\max_Vu> -A=\log\min\le\{1,\f{1}{4|V|\Lambda^2}\ri\}.$$
This ends the proof of the lemma, and completes the proof of Theorem \ref{prior-1}. $\hfill\Box$

\section{Topological degree}\label{sec-Brouwer}

In this section, we shall prove Theorem \ref{degree}. Precisely we shall compute the topological degree of certain maps related to
the Chern-Simons Higgs model.\\

{\it Proof of Theorem \ref{degree}.}
Assume $V=\{x_1,\cdots,x_\ell\}$. let $X=L^\infty(V)$. We may identify $X$ with the Euclidean space $\mathbb{R}^\ell$.
Without causing ambiguity, we define a map $F:X\times[0,1]\ra X$ by
$$F(u,\sigma)=-\Delta u+\lambda e^u(e^u-\sigma)+f,\quad (u,\sigma)\in X\times [0,1].$$
Obviously, $F$ is a smooth map. For the fixed real number $\lambda$ and the fixed function $f$, since $\lambda\overline{f}\not=0$,
 there must exist a large number $\Lambda>0$ such that
\be\label{bound}\Lambda^{-1}\leq |\lambda|\leq \Lambda,\,\,\,
\Lambda^{-1}\leq \le|\int_Vfd\mu\ri|\leq \Lambda,\,\,\,\|f\|_{L^\infty(V)}\leq \Lambda.\ee
Here and in the sequel, $\overline{f}$ denotes the integral mean of a function $f$.
Then it follows from Theorem \ref{prior-1} that
there exists a constant $R_0>0$, depending only on $\Lambda$ and the graph $V$, such that for all $\sigma\in[0,1]$, all solutions of $F(u,\sigma)=0$
satisfy $\|u\|_{L^\infty(V)}< R_0$. Denote a ball centered at $0\in X$ with radius $r$ by $B_r\subset X$, and its boundary by
$\p B_r=\{u\in X:\|u\|_{L^\infty(V)}=r\}$. Thus we conclude
$$0\not\in F(\p B_{R},\sigma),\quad\forall \sigma\in[0,1],\,\,\forall R\geq R_0.$$
By the homotopic invariance of the topological degree, we have
\be\label{degree-1}\deg(F(\cdot,1),B_{R},0)=\deg(F(\cdot,0),B_{R},0),\quad\forall R\geq R_0.\ee
Given any $\epsilon>0$, we define another smooth map $G_\epsilon:X\times[0,1]\ra X$ by
$$G_\epsilon(u,t)=-\Delta u+\lambda e^{2u}+(t+(1-t)\epsilon) f,\quad (u,t)\in X\times [0,1].$$
Notice that
$$\le|(t+(1-t)\epsilon)\int_Vfd\mu\ri|\geq\min\{1,\epsilon\}\le|\int_Vfd\mu\ri|,\quad\forall t\in[0,1].$$
Applying Theorem \ref{prior-1} again, we find a constant $R_\epsilon>0$, depending only on $\epsilon$, $\Lambda$ and the graph
$V$, such that
all solutions $u$ of $G_\epsilon(u,t)=0$ satisfy $\|u\|_{L^\infty(V)}< R_\epsilon$ for all $t\in [0,1]$. This implies
$$0\not\in G_\epsilon(\p B_{R_\epsilon},t),\quad\forall t\in[0,1].$$
Hence the homotopic invariance of the topological degree leads to
\be\label{degree-2}\deg(G_\epsilon(\cdot,1),B_{R_\epsilon},0)=\deg(G_\epsilon(\cdot,0),B_{R_\epsilon},0).\ee
To calculate $\deg(G_\epsilon(\cdot,0),B_{R_\epsilon},0)$, we need to understand the solvability of the equation
\be\label{eq-4}G_\epsilon(u,0)=-\Delta u+\lambda e^{2u}+\epsilon f=0.\ee
 Now we {\it claim} two properties of solutions of (\ref{eq-4}):
 $(i)$ If $\lambda\overline{f}<0$, then there exists an $\epsilon_0>0$ such that for any $\epsilon\in (0,\epsilon_0)$,
(\ref{eq-4}) has a unique solution $u_\epsilon$, which satisfies $e^{2u_\epsilon}\leq C\epsilon$, where $C$ is a constant depending only on
$\Lambda$ and the graph $V$; $(ii)$ If $\lambda\overline f>0$, then (\ref{eq-4}) has no solution for all $\epsilon>0$.

To see Claim $(i)$,
 for any $\epsilon>0$, we let $v_\epsilon$ be the unique solution of the equation
$$\le\{\begin{array}{lll}
\Delta v=\epsilon f-\epsilon \overline{f}\quad{\rm in}\quad V\\[1.2ex]
\overline{v}=0.
\end{array}\ri.$$
Then the solvability of (\ref{eq-4}) is equivalent to that of the equation
\be\label{20}\Delta w=\lambda e^{2v_\epsilon}e^{2w}+\epsilon\overline{f}.\ee
Note that the existence of solutions to (\ref{20}), under the assumptions that $\epsilon$ is sufficiently small and  $\lambda\overline{f}<0$, follows from (\cite{Gri1}, Theorems 2 and 4). Hence there exists some $\epsilon_1>0$ such that if $0<\epsilon<\epsilon_1$, then the equation (\ref{eq-4}) has a solution $u_\epsilon$.
 Integrating both sides of (\ref{eq-4}), we have by (\ref{bound}),
$$\int_Ve^{2u_\epsilon}d\mu=-\f{\epsilon}{\lambda}\int_Vfd\mu\leq \Lambda^2\epsilon,$$
which leads to
\be\label{small}e^{2u_\epsilon(x)}\leq \f{\Lambda^2}{\mu_0}\epsilon,\quad\forall x\in V,\ee
where $\mu_0=\min_{x\in V}\mu(x)$. We also need to prove the uniqueness of the solution. Let $\varphi$ be an arbitrary solution of (\ref{eq-4}), namely it satisfies
\be\label{var-eq}\Delta\varphi=\lambda e^{2\varphi}+\epsilon f.\ee
The same procedure as above gives
\be\label{var}\int_Ve^{2\varphi}d\mu\leq \Lambda^2\epsilon,\quad
e^{2\varphi(x)}\leq \f{\Lambda^2}{\mu_0}\epsilon\quad{\rm for\,\,\,all}\quad x\in V.\ee
Subtracting (\ref{var-eq}) from (\ref{eq-4}) and integrating by parts, we have
$$0=\int_V\Delta(u_\epsilon-\varphi)d\mu=\lambda\int_V(e^{2u_\epsilon}-e^{2\varphi})d\mu,$$
which leads to
$$\min_V(u_\epsilon-\varphi)<0<\max_V(u_\epsilon-\varphi).$$
As a consequence, there holds
\be\label{abs}|u_\epsilon-\varphi|\leq \max_V(u_\epsilon-\varphi)-\min_V(u_\epsilon-\varphi).\ee
Also we derive from (\ref{eq-4}), (\ref{small}), (\ref{var-eq}), and (\ref{var}),
\bea\nonumber
|\Delta(u_\epsilon-\varphi)(x)|&=& \le|\lambda \le(e^{2u_\epsilon(x)}- e^{2\varphi(x)}\ri)\ri|\\\nonumber
&\leq&2\Lambda\le(e^{2u_\epsilon(x)}+ e^{2\varphi(x)}\ri)|u_\epsilon(x)-\varphi(x)|\\\label{est}
&\leq&\f{4\Lambda^3}{\mu_0}\epsilon|u_\epsilon(x)-\varphi(x)|.
\eea
Combining (\ref{equiv}), (\ref{abs}) and (\ref{est}), we obtain
\be\label{contr}\max_V(u_\epsilon-\varphi)-\min_V(u_\epsilon-\varphi)\leq \sqrt{\f{(\ell-1)|V|}{w_0\lambda_1}}\f{4\Lambda^3}{\mu_0}\epsilon
\le(\max_V(u_\epsilon-\varphi)-\min_V(u_\epsilon-\varphi)\ri).\ee
Choose
$$\epsilon_0=\min\le\{\epsilon_1,\sqrt{\f{w_0\lambda_1}{(\ell-1)|V|}}\f{\mu_0}{8\Lambda^3}\ri\}.$$
If we take $0<\epsilon<\epsilon_0$, then (\ref{contr}) implies $\varphi\equiv u_\epsilon$ on $V$, and thus (\ref{eq-4}) has a unique solution. Hence $(i)$ holds.

To see Claim $(ii)$, in the case $\lambda\overline{f}>0$, if (\ref{eq-4}) has a solution $u$, then there holds
$$0=\int_V\Delta ud\mu=\lambda\int_Ve^{2u}d\mu+\int_Vfd\mu,$$
which is impossible. This confirms $(ii)$, and our claims hold.

Let us continue to prove the theorem.
Note that $-\Delta: X\ra X$ is a nonnegative definite symmetric operator, its eigenvalues are written as
$$0=\lambda_0<\lambda_1\leq \lambda_2\leq\cdots\leq \lambda_{\ell-1},$$
where $\ell$ is the number of all points in $V$.
By Claim $(i)$, in the case $\lambda\overline{f}<0$, we may choose a sufficiently small $\epsilon>0$ such that $G_\epsilon(u,0)=0$ has a unique
solution $u_\epsilon$ verifying
$$2|\lambda|e^{2u_\epsilon(x)}<\lambda_1.$$
A straightforward calculation shows
$$DG_\epsilon(u_\epsilon,0)=-\Delta+2\lambda e^{2u_\epsilon}{\rm I},$$
where we identify the linear operator $-\Delta$ with the $\ell\times\ell$ matrix corresponding to $-\Delta$, and denote
the $\ell\times\ell$ diagonal matrix ${\rm diag}[1,1,\cdots,1]$ by ${\rm I}$.
Clearly
$$\deg(G_\epsilon(\cdot,0),B_{R_\epsilon},0)={\rm sgn\,det}\le(DG_\epsilon(u_\epsilon,0)\ri)
={\rm sgn}\le\{2\lambda e^{2u_\epsilon(x)}\Pi_{j=1}^{\ell-1}(\lambda_j+2\lambda e^{2u_\epsilon(x)})\ri\}={\rm sgn}\lambda.$$
This together with (\ref{degree-1}) and (\ref{degree-2}) leads to
\bna
\deg(F(\cdot,1),B_{R_\epsilon},0)&=&\deg(F(\cdot,0),B_{R_\epsilon},0)\\
&=&\deg(G_\epsilon(\cdot,1),B_{R_\epsilon},0)\\
&=&\deg(G_\epsilon(\cdot,0),B_{R_\epsilon},0)\\
&=&{\rm sgn}\lambda.
\ena
By Claim $(ii)$, in the case $\lambda\overline{f}>0$, since $G_\epsilon(u,0)=0$ has no solution, we obtain
$$\deg(F(\cdot,1),B_{R_\epsilon},0)=\deg(G_\epsilon(\cdot,0),B_{R_\epsilon},0)=0.$$
Thus the proof of Theorem \ref{degree} is completed. $\hfill\Box$

\section{Existence results}\label{sec-existence}

In this section, we shall prove Theorem \ref{existence} by using the topological degree in Theorem \ref{degree}.\\

{\it Proof of Theorem \ref{existence} $(\mathsf{a})$.}  If $\lambda\overline{f}<0$, then by Theorem \ref{degree}, we find some large $R_0>1$ such that
$$\deg(F,B_{R_0},0)\not=0.$$
Thus the Kronecker's existence theorem implies (\ref{eq-1}) has a solution. $\hfill\Box$\\

In the remaining part of this section, we always assume $\lambda\overline{f}>0$. We first prove that
(\ref{eq-1}) has a local minimum solution for large $|\lambda|$, say

\begin{lemma}\label{large}
If $|\lambda|$ is chosen sufficiently large, then the equation $(3)$ has a local minimum solution.
\end{lemma}
\proof Let us first consider the subcase $\lambda>0$ and $\overline{f}>0$.
Set
\be\label{Llambda}L_\lambda u=-\Delta u+\lambda e^u(e^u-1)+f.\ee
For real numbers $A$ and $\lambda$, there hold
$$L_\lambda A=\lambda e^A(e^A-1)+f,\quad
L_\lambda \log \f{1}{2}=-\f{1}{4}\lambda +f.$$
Clearly, taking sufficiently large $A>1$ and $\lambda>1$, we have
\be\label{u-l}L_\lambda A>0,\quad L_\lambda \log\f{1}{2}<0.\ee
Recall the functional $J_\lambda:X=L^\infty(V)\ra\mathbb{R}$ defined by (\ref{funct-1}). Since $X\cong \mathbb{R}^\ell$,
$J_\lambda\in C^2(X,\mathbb{R})$, and $\{u\in X:\log\f{1}{2}\leq u\leq A\}$ is a bounded closed subset of $X$, it is easy to find some $u_\lambda\in X$
satisfying $\log \f{1}{2}\leq u_\lambda(x)\leq A$ for all $x\in V$ and
\be\label{J-min}J_\lambda(u_\lambda)=\min_{\log\f{1}{2}\leq u\leq A}J_\lambda(u).\ee
We {\it claim} that
\be\label{strict}\log\f{1}{2}< u_\lambda(x)< A\quad{\rm for\,\,all}\quad x\in V.\ee
Suppose not. There must hold $u_\lambda(x_0)=\log\f{1}{2}$ for some $x_0\in V$, or $u_\lambda(x_1)=A$ for some $x_1\in V$.
If $u_\lambda(x_0)=\log\f{1}{2}$, then we take a small $\epsilon>0$ such that
$$\log\f{1}{2}<u_\lambda(x)+t\delta_{x_0}(x)<A,\quad\forall x\in V,\,\forall t\in(0,\epsilon).$$
On one hand, in view of (\ref{u-l}) and (\ref{J-min}), we have
\bea\nonumber
0&\leq&\le.\f{d}{dt}\ri|_{t=0}J_\lambda(u_\lambda+t\delta_{x_0})\\\nonumber
&=&\int_V\le(-\Delta u_\lambda+\lambda e^{u_\lambda}(e^{u_\lambda}-1)+f\ri)\delta_{x_0}d\mu\\\nonumber
&=&-\Delta u_\lambda(x_0)+\lambda e^{u_\lambda(x_0)}(e^{u_\lambda(x_0)}-1)+f(x_0)\\\label{less}
&<&-\Delta u_\lambda(x_0).
\eea
On the other hand, since $u_\lambda(x)\geq u_\lambda(x_0)$ for all $x\in V$, we conclude $\Delta u_\lambda(x_0)\geq 0$,
which contradicts (\ref{less}). Hence $u_\lambda(x)>\log\f{1}{2}$ for all $x\in V$. In the same way, we exclude the possibility of
$u_\lambda(x_1)=A$ for some $x_1\in V$. This confirms our claim (\ref{strict}). Combining (\ref{J-min}) and (\ref{strict}), we conclude
that $u_\lambda$ is a local minimum critical point of $J_\lambda$, in particular, $u_\lambda$ is a solution of (\ref{eq-1}).

Now we consider the subcase $\lambda<0$ and $\overline{f}<0$. Let $\varphi$ be the unique solution of
$$\le\{\begin{array}{lll}
\Delta\varphi=f-\overline{f}\\[1.2ex]
\overline{\varphi}=0.
\end{array}\ri.$$
Using the notation of the operator $L_\lambda$ given by (\ref{Llambda}), we have
\bea\nonumber
L_\lambda(\varphi-A)&=&-\Delta\varphi+\lambda e^{\varphi-A}(e^{\varphi-A}-1)+f\\\nonumber
&=&\lambda e^{\varphi-A}(e^{\varphi-A}-1)+\overline{f}\\\label{lo-bd}
&<&0
\eea
and
\bna
L_\lambda(\log\f{1}{2})&=&\lambda e^{\log\f{1}{2}}(e^{\log\f{1}{2}}-1)+f\\
&=&-\f{\lambda}{4}+f\\
&>&0,
\ena
provided that $\lambda<4\min_Vf$ and $A>1$ is chosen sufficiently large. Similar to (\ref{J-min}) and (\ref{strict}), there exists
some $u_\lambda$ satisfying $\varphi(x)-A< u_\lambda(x)< \log\f{1}{2}$ for all $x\in V$ and
$$J_\lambda(u_\lambda)=\min_{\varphi-A\leq u\leq \log\f{1}{2}}J_\lambda(u)=\min_{\varphi-A< u< \log\f{1}{2}}J_\lambda(u).$$
This implies $u_\lambda$ is a local minimum solution of (\ref{eq-1}). $\hfill\Box$\\

To proceed, we also need the following:

\begin{lemma}\label{solution}
If $\lambda_1>0$ such that the equation $L_{\lambda_1}u=0$ has a solution $u_{\lambda_1}$, then for any $\lambda>\lambda_1$, we have
$$L_{\lambda} \le(u_{\lambda_1}+\log\f{\lambda_1}{\lambda}\ri)<0.$$
Similarly, if $\lambda_2<0$ such that $L_{\lambda_2}u_{\lambda_2}=0$, then for any $\lambda<\lambda_2$, there hods
$$L_{\lambda} \le(u_{\lambda_2}+\log\f{\lambda_2}{\lambda}\ri)>0.$$
\end{lemma}
\proof If $\lambda>\lambda_1>0$, then
\bna
L_\lambda \le(u_{\lambda_1}+\log\f{\lambda_1}{\lambda}\ri)&=&-\Delta u_{\lambda_1}+\lambda_1 e^{u_{\lambda_1}}
\le(\f{\lambda_1}{\lambda}e^{u_{\lambda_1}}-1\ri)+f \\
&<&-\Delta u_{\lambda_1}+\lambda_1 e^{u_{\lambda_1}}
(e^{u_{\lambda_1}}-1)+f\\
&=&0.
\ena
If $\lambda<\lambda_2<0$, then
\bna
L_\lambda \le(u_{\lambda_2}+\log\f{\lambda_2}{\lambda}\ri)&=&-\Delta u_{\lambda_2}+\lambda_2 e^{u_{\lambda_2}}
\le(\f{\lambda_2}{\lambda}e^{u_{\lambda_2}}-1\ri)+f\\
&>&-\Delta u_{\lambda_2}+\lambda_2 e^{u_{\lambda_2}}
\le(e^{u_{\lambda_2}}-1\ri)+f\\
&=&0,
\ena
as we desired. $\hfill\Box$\\

As a consequence, we have

\begin{lemma}\label{minimum-sol}
Assume $L_{\lambda_1}u_{\lambda_1}=L_{\lambda_2}u_{\lambda_2}=0$ on $V$. If either $\lambda>\lambda_1>0$ or $\lambda<\lambda_2<0$,
then the equation (3) has a local minimum solution $u_\lambda$.
\end{lemma}
\proof Assume $\lambda>\lambda_1>0$. Let $A>1$ be a sufficiently large constant such that $L_\lambda A>0$ and $u_{\lambda_1}+\log\f{\lambda_1}{\lambda}< A$
on $V$. Then there exists some $u_\lambda$ such that
$$J_\lambda(u_\lambda)=\min_{u_{\lambda_1}+\log\f{\lambda_1}{\lambda}\leq u\leq A}J_\lambda(u).$$
Suppose there is some point $x_0\in V$ satisfying $u_\lambda(x_0)=u_{\lambda_1}(x_0)+\log\f{\lambda_1}{\lambda}$.
Let $\epsilon>0$ be so small that for $t\in(0,\epsilon)$, there holds
$$u_{\lambda_1}(x)+\log\f{\lambda_1}{\lambda}< u_\lambda(x)+t\delta_{x_0}(x)< A\quad{\rm for\,\,all}\quad x\in V.$$
 Similarly as
we did in the proof of Lemma \ref{large}, we have by Lemma \ref{solution},
\bna
0&\leq&\le.\f{d}{dt}\ri|_{t=0}J_\lambda(u_\lambda+t\delta_{x_0})\\
&=&-\Delta u_\lambda(x_0)+\lambda e^{u_\lambda(x_0)}(e^{u_\lambda(x_0)}-1)+f(x_0)\\
&=&-\Delta\le(u_\lambda-u_{\lambda_1}\ri)(x_0)+L_{\lambda}\le(u_{\lambda_1}+\log{\f{\lambda_1}{\lambda}}\ri)(x_0)\\
&<&-\Delta\le(u_\lambda-u_{\lambda_1}\ri)(x_0).
\ena
This contradicts the fact that $x_0$ is a minimum point of $u_\lambda-u_{\lambda_1}-\log\f{\lambda_1}{\lambda}$. Hence
$$u_\lambda(x)>u_{\lambda_1}(x)+\log\f{\lambda_1}{\lambda},\quad\forall x\in V.$$
In the same way we obtain $u(x)<A$ for all $x\in V$. Therefore $u_\lambda$ is a local minimum critical point of $J_\lambda$.

Assume $\lambda<\lambda_2<0$. The constant $A>1$ is chosen sufficiently large such that
$\varphi-A<u_{\lambda_2}+\log\f{\lambda_2}{\lambda}$ on $V$, and  $\varphi-A$ satisfies (\ref{lo-bd}). Clearly there exists some
$u_\lambda$ such that
$$J_\lambda(u_\lambda)=\min_{\varphi-A\leq u\leq u_{\lambda_2}+\log\f{\lambda_2}{\lambda}}J_\lambda(u).$$
If there is some point $x_1\in V$ satisfying $u_\lambda(x_1)=u_{\lambda_2}(x_1)+\log\f{\lambda_2}{\lambda}$, then there is a small
$\epsilon>0$ such that for $t\in(0,\epsilon)$, there holds
$$\varphi(x)-A< u_\lambda(x)-t\delta_{x_1}(x)< u_{\lambda_2}(x)+\log\f{\lambda_2}{\lambda}\quad{\rm for\,\,all}\quad x\in V.$$
 Thus we have
by Lemma \ref{solution},
\bna
0&\leq&\le.\f{d}{dt}\ri|_{t=0}J_\lambda(u_\lambda-t\delta_{x_1})\\
&=&\Delta u_\lambda(x_1)-\lambda e^{u_\lambda(x_1)}(e^{u_\lambda(x_1)}-1)-f(x_0)\\
&=&\Delta\le(u_\lambda-u_{\lambda_2}\ri)(x_0)-L_{\lambda}\le(u_{\lambda_2}+\log{\f{\lambda_2}{\lambda}}\ri)(x_0)\\
&<&\Delta\le(u_\lambda-u_{\lambda_2}\ri)(x_0).
\ena
This contradicts the fact that $x_1$ is a maximum point of $u_\lambda-u_{\lambda_2}-\log\f{\lambda_2}{\lambda}$. Hence
$$u_\lambda(x)<u_{\lambda_2}(x)+\log\f{\lambda_2}{\lambda},\quad\forall x\in V.$$
In the same way we obtain $u(x)>\varphi(x)-A$ for all $x\in V$. Therefore $u_\lambda$ is a local minimum critical point of $J_\lambda$.
Thus we complete the proof of the lemma. $\hfill\Box$\\

We conclude from Lemmas \ref{large} and \ref{minimum-sol} that the following two critical numbers are well defined.
\bea
\label{g-1}
&\Lambda^\ast=\inf\le\{\lambda>0: \lambda\overline{f}>0, J_\lambda \,\,{\rm has\,\,a\,\,local\,\,minimum\,\,critical\,\,point}\ri\}\\[1.2ex]
&\Lambda_\ast=\sup\le\{\lambda<0: \lambda\overline{f}>0, J_\lambda \,\,{\rm has\,\,a\,\,local\,\,minimum\,\,critical\,\,point}\ri\}.\label{g-2}
\eea
\begin{lemma}\label{lambda}
If $\overline{f}>0$, then $\Lambda^\ast\geq 4\overline{f}$; If $\overline{f}<0$, then $\Lambda_\ast\leq 4\overline{f}$.
\end{lemma}
\proof Suppose $\lambda\not=0$ and $u$ is a solution of
$\Delta u=\lambda e^u(e^u-1)+f$. Integration by parts gives
$$-\f{\int_Vfd\mu}{\lambda}=\int_Ve^u(e^u-1)d\mu\geq -\f{|V|}{4},$$
since $e^u(e^u-1)\geq -\f{1}{4}$. The conclusion follows from (\ref{g-1}) and (\ref{g-2}) immediately. $\hfill\Box$\\

We are now ready to complete the proof of the remaining part of the theorem. \\

{\it Proof of Theorem \ref{existence} $(\mathsf{b})$.} {\it We first consider the solvability of the equation (\ref{eq-1}) under the assumption $\lambda\in (0,\Lambda^\ast]\cup[\Lambda_\ast,0)$.}

If $\lambda\in (0,\Lambda^\ast)\cup(\Lambda_\ast,0)$, then $(\ref{eq-1})$ has no solution. Indeed,
suppose there exists a number $\lambda_1\in (0,\Lambda^\ast)\cup(\Lambda_\ast,0)$ such that
(\ref{eq-1}) has a solution at $\lambda=\lambda_1$. With no loss of generality, we assume $\lambda_1\in (\Lambda_\ast,0)$, then by Lemma \ref{minimum-sol},
(\ref{eq-1}) has a local minimum solution at any $\lambda\in [\Lambda_\ast,\lambda_1)$. This contradicts the definition of $\Lambda_\ast$.
Hence $(\ref{eq-1})$ has no solution for any $\lambda\in (0,\Lambda^\ast)\cup(\Lambda_\ast,0)$.

Note that for any $j\in\mathbb{N}$, there exists a solution $u_j$ of (\ref{eq-1}) with $\lambda=\Lambda_\ast-1/j$.
According to Theorem \ref{prior-1}, $(u_j)$ is uniformly bounded in $V$. Thus up to a subsequence, $(u_j)$ uniformly
converges to some function $u^\ast$, a solution of (\ref{eq-1}) with $\lambda=\Lambda_\ast$. In the same way,
(\ref{eq-1}) has also a solution at $\lambda=\Lambda^\ast$.\\

{\it We next consider multiple solutions of (\ref{eq-1}) under the assumption $\lambda\in(\Lambda^\ast,+\infty)\cup(-\infty,\Lambda_\ast)$.}

If $\lambda \in(\Lambda^\ast,+\infty)\cup(-\infty,\Lambda_\ast)$, by (\ref{g-1}) and (\ref{g-2}),
we let $u_\lambda$ be a local minimum critical point of $J_\lambda$. With no loss of generality, we may assume $u_\lambda$ is the unique
critical point of $J_\lambda$. For otherwise, $J_\lambda$ has already at least two critical points, and the proof terminates.
According to (\cite{Chang93}, Chapter 1, page 32), the $q$-th critical group of ${J}_\lambda$ at $u_\lambda$ is defined by
 \be\label{group}\mathsf{C}_q({J}_\lambda,u_\lambda)=\mathsf{H}_q({J}_\lambda^c\cap {U},\{{J}_\lambda^c\setminus\{u_\lambda\}\}
 \cap {U},\mathsf{G}),\ee
 where ${J}_\lambda(u_\lambda)=c$, ${J}_\lambda^c=\{u\in X:{J}_\lambda(u)\leq c\}$, $U$ is a neighborhood of
 $u_\lambda\in X$, $\mathsf{H}_q$
 is the singular homology group with the coefficients groups $\mathsf{G}$, say $\mathbb{Z}$, $\mathbb{R}$. By the excision property of
 $\mathsf{H}_q$, this definition is not dependent on the choice of $U$. It is easy to calculate
\be\label{cri}{\mathsf{C}}_q(J_\lambda,u_\lambda)=\delta_{q0}\mathsf{G}.\ee
Note that $J_\lambda$ satisfies the Palais-Smale condition. Indeed, if $J_\lambda(u_j)\ra c\in\mathbb{R}$ and $J^\prime(u_j)\ra 0$
as $j\ra\infty$, then using the method of proving Theorem \ref{prior-1}, we obtain $(u_j)$ is uniformly bounded. Since $X$ is pre-compact,
then up to a subsequence, $(u_j)$ converges uniformly to some $u^\ast$, a critical point $J_\lambda$. Thus the Palais-Smale condition follows.
Notice also that
$$DJ_\lambda(u)=-\Delta u+\lambda e^{u}(e^u-1)+f=F(u),$$
where $F$ is given as in Theorem \ref{degree}.
According to (\cite{Chang93}, Chapter 2, Theorem 3.2), in view of (\ref{cri}), we have for sufficiently large $R>1$,
$$\deg(F,B_{R},0)=\deg(DJ_\lambda,B_R,0)=\sum_{q=0}^\infty (-1)^q{\rm rank}\, \mathsf{C}_q(J_\lambda,u_{\lambda})=1.$$
This contradicts $\deg(F,B_{R},0)=0$ derived from Theorem \ref{degree}. Therefore the equation (\ref{eq-1}) has at least two different solutions,
and the proof of Theorem \ref{existence} $(\mathsf{b})$ is finished. $\hfill\Box$

\section{Chern-Simons Higgs System}\label{sec-system}
In this section, we shall calculate the topological degree of the map related to the Chern-Simons Higgs system (\ref{system}), and then
use the degree to obtain partial results for multiplicity of solutions to the system. In particular,
Theorems \ref{system-apriori}-\ref{syst-thm} will be proved. We first derive a priori estimate for solutions of (\ref{system-8}), a deformation of (\ref{system}). \\

{\it Proof of Theorem \ref{system-apriori}.}
Let $\sigma\in[0,1]$, $\lambda>0$, $\eta>0$, $\overline{f}>0$, $\overline{g}>0$, and $(u,v)$ be a solution of the system
(\ref{system-8}). Note that there exist a unique solution $\varphi$ to the equation
$$\le\{
\begin{array}{lll}
\Delta \varphi=f-\overline{f}\\[1.2ex]
\int_V\varphi d\mu=0
\end{array}\ri.
$$
and a unique solution $\psi$ to the equation
$$\le\{
\begin{array}{lll}
\Delta \psi=g-\overline{g}\\[1.2ex]
\int_V\psi d\mu=0.
\end{array}\ri.
$$
Set $w=u-\varphi$ and $z=v-\psi$. Then we have
\be\label{system-2}
\le\{
\begin{array}{lll}
\Delta w=\lambda e^\psi e^z(e^{\varphi}e^w-\sigma)+\overline{f}\\[1.2ex]
\Delta z=\eta e^{\varphi}e^w(e^\psi e^z-\sigma)+\overline{g},
\end{array}\ri.
\ee
We {\it claim} that
\be\label{up-2}w(x)<-\min_V\varphi\quad {\rm for\,\, all}\quad x\in V.\ee Suppose not. There necessarily hold
$\max_Vw\geq -\min_V\varphi$.
Take $x_0\in V$ satisfying $w(x_0)=\max_Vw$. Since $\sigma\in[0,1]$, $\lambda>0$, $\overline{f}>0$
and $\varphi(x_0)+w(x_0)\geq 0$,
we have
$$0\geq \Delta w(x_0)=\lambda e^{\psi(x_0)} e^{z(x_0)}(e^{\varphi(x_0)}e^{w(x_0)}-\sigma)+\overline{f}\geq\overline{f}> 0,$$
which is impossible. Hence our claim (\ref{up-2}) follows.  Keeping in mind $\eta>0$ and $\overline{g}>0$, in the same way as above,
we also have
\be\label{up-4}z(x)<-\min_V\psi\quad {\rm for\,\, all}\quad x\in V.\ee
Inserting (\ref{up-2}) and (\ref{up-4}) into (\ref{system-2}), we obtain
$$\|\Delta w\|_{L^\infty(V)}+\|\Delta z\|_{L^\infty(V)}\leq C$$
for some constant $C$, depending only on $\lambda,\eta,f,g$ and the graph $V$. The most important thing here is that the
constant $C$
is not dependent on the parameter $\sigma\in[0,1]$. Coming back to the inequality (\ref{equiv}), we immediately conclude
\be\label{diffe-1}\max_V w-\min_V w\leq C\ee
and
$$\max_V z-\min_V z\leq C.$$
Observe that integration on both sides of the second equation in (\ref{system-2}) leads to
$$\int_Ve^{\varphi}e^w(e^\psi e^z-\sigma)d\mu=-\f{\overline{g}}{\eta}|V|.$$
As a consequence, there holds
\bna
0<\f{\overline{g}}{\eta}\leq e^{\max_Vw}e^{\max_V\varphi}\le(e^{\max_V\psi}+1\ri)
\leq Ce^{\max_Vw}.
\ena
Hence $\max_Vw\geq -C$,
and in view of (\ref{diffe-1}),
\be\label{lower}\min_Vw\geq -C.\ee
In the same way, from (\ref{diffe-1}) and the first equation of (\ref{system-2}), we derive
\be\label{lower-2}\min_Vz\geq -C.\ee
In view of (\ref{up-2}), (\ref{up-4}), (\ref{lower}) and (\ref{lower-2}), the proof of the theorem
is completed. $\hfill\Box$ \\

Now we calculate the topological degree of the map defined as in (\ref{map-syst}).\\

{\it Proof of Theorem \ref{degree-system}.}
Let $X=L^\infty(V)$. Define a map $\mathcal{F}:X\times X\times[0,1]\ra X\times X$ by
$$\mathcal{F}(u,v,\sigma)=(-\Delta u+\lambda e^v(e^u-\sigma)+f,-\Delta v+\eta e^u(e^v-\sigma)+g),\quad\forall
(u,v,\sigma)\in X\times X\times [0,1].$$
Obviously $\mathcal{F}\in C^2(X\times X\times[0,1],X\times X)$.
On one hand, by Theorem \ref{system-apriori}, there exists some $R_0>0$ such that for any $R\geq R_0$, we have
$$0\not\in \mathcal{F}(\p B_R,\sigma),\quad\forall \sigma\in[0,1],$$
and thus the homotopic invariance of the topological degree implies
\be\label{homotopy}\deg(\mathcal{F}(\cdot,1),B_R,(0,0))=\deg(\mathcal{F}(\cdot,0),B_R,(0,0)).\ee
Here we denote $B_R=\{(u,v)\in X\times X:\|u\|_{L^\infty(V)}+\|v\|_{L^\infty(V)}<R\}$ and
$\p B_R=\{(u,v)\in X\times X:\|u\|_{L^\infty(V)}+\|v\|_{L^\infty(V)}=R\}$, as usual.

On the other hand, we calculate $\deg(\mathcal{F}(\cdot,0),B_R,(0,0))$.
Since $\lambda>0$ and $\overline{f}>0$, integrating both sides
of the first equation of the system
\be\label{system-0}
\le\{
\begin{array}{lll}
\Delta u=\lambda e^{u+v}+f\\[1.2ex]
\Delta v=\eta e^{u+v}+g,
\end{array}\ri.
\ee
we get a contradiction, provided that (\ref{system-0}) is solvable. This implies
$$\le\{(u,v)\in X\times X: \mathcal{F}(u,v,0)=(0,0)\ri\}=\varnothing.$$
As a consequence, there holds
\be\label{0-deg}\deg(\mathcal{F}(\cdot,0),B_R,(0,0))=0.\ee
Combining (\ref{homotopy}) and (\ref{0-deg}), we get the desired result. $\hfill\Box$\\

Let $\mathcal{J}_\lambda:X\times X\ra \mathbb{R}$ be a functional defined as in (\ref{functional}).
Note that the critical point of $\mathcal{J}_\lambda$ is a solution of the Chern-Simons system (\ref{system}).
The following property of $\mathcal{J}_\lambda$ will be not only useful for our subsequent analysis, but also
of its own interest.

\begin{lemma}\label{PS}
Under the assumptions $\lambda>0$, $\overline{f}>0$ and $\overline{g}>0$,
$\mathcal{J}_\lambda$ satisfies the Palais-Smale condition at any level $c\in\mathbb{R}$.
\end{lemma}

\proof Let $c\in\mathbb{R}$ and $\{(u_k,v_k)\}$ be a sequence in $X\times X$ such that $\mathcal{J}_\lambda(u_k,v_k)\ra c$ and
$$\mathcal{J}_\lambda^\prime(u_k,v_k)\ra (0,0)\quad{\rm in}\quad (X\times X)^\ast\cong \mathbb{R}^{\ell}\times\mathbb{R}^\ell.$$
This together with  (\ref{derivative}) gives
\be\label{deriv}\le\{\begin{array}{lll}
-\Delta u_k+\lambda e^{v_k}(e^{u_k}-1)+f=o_k(1)\\[1.2ex]
-\Delta v_k+\lambda e^{u_k}(e^{v_k}-1)+g=o_k(1),
\end{array}\ri.\ee
where $o_k(1)\ra 0$ uniformly on $V$ as $k\ra \infty$. Comparing (\ref{deriv}) with the system (\ref{system}),
we have by using the same method
as in the proof of Theorem \ref{system-apriori},
$$\|u_k\|_{L^\infty(V)}+\|v_k\|_{L^\infty(V)}\leq C$$
for some constant $C$, provided that $k\geq k_1$ for some large positive integer $k_1$. Since $V$ is finite, $X$ is pre-compact.
Hence, up to a subsequence, $u_k\ra u^\ast$ and $v_k\ra v^\ast$ uniformly in $V$ for some functions $u^\ast$ and $v^\ast$.
Obviously $\mathcal{J}_\lambda^\prime(u^\ast,v^\ast)=(0,0)$. Thus $\mathcal{J}_\lambda$ satisfies the $(PS)_c$ condition.
$\hfill\Box$ \\

Finally we prove a partial multiple solutions result for the system (\ref{system}).\\

{\it Proof of Theorem \ref{syst-thm}}. We distinguish two hypotheses to proceed.

{\it Case} 1. {\it $\mathcal{J}_\lambda$ has a non-degenerate critical point $(u_\lambda,v_\lambda)$.}

Since $(u_\lambda,v_\lambda)$ is non-degenerate, we have
$$\det D^2\mathcal{J}_\lambda(u_\lambda,v_\lambda)\not =0.$$
 Suppose $(u_\lambda,v_\lambda)$ is the unique critical point of $\mathcal{J}_\lambda$.
 Then we conclude for all $R>\|u_\lambda\|_{L^\infty(V)}+\|v_\lambda\|_{L^\infty(V)}$,
 \be\label{non-0}\deg(D\mathcal{J}_\lambda, B_R,(0,0))={\rm sgn}\,\det  D^2\mathcal{J}_\lambda(u_\lambda,v_\lambda)\not=0.\ee
 Here and in the sequel, as in the proof of Theorem \ref{degree-system}, $B_R$ is a ball centered at $(0,0)$ with radius $R$.
   Notice that $D\mathcal{J}_\lambda(u,v)=\mathcal{F}(u,v)$ for all $(u,v)\in X\times X$, where $\mathcal{F}$ is defined as in
 (\ref{map-syst}). By Theorem \ref{degree-system}, we have
 $$\deg(D\mathcal{J}_\lambda,B_R,(0,0))=\deg(\mathcal{F},B_R,(0,0))=0,$$
 contradicting (\ref{non-0}). Hence $\mathcal{J}_\lambda$ must have at least two critical points.

 {\it Case $2$. $\mathcal{J}_\lambda$ has a local minimum critical point $(\varphi_\lambda,\psi_\lambda)$.}

 Similar to (\ref{group}), the $q$-th critical group of $\mathcal{J}_\lambda$ at the critical point $(\varphi_\lambda,\psi_\lambda)$ reads as
 $$\mathsf{C}_q(\mathcal{J}_\lambda,(\varphi_\lambda,\psi_\lambda))=\mathsf{H}_q(\mathcal{J}_\lambda^c\cap \mathscr{U},\{\mathcal{J}_\lambda^c\setminus\{(\varphi_\lambda,\psi_\lambda)\}\}
 \cap\mathscr{U},\mathsf{G}),$$
 where $\mathcal{J}_\lambda(\varphi_\lambda,\psi_\lambda)=c$, $\mathcal{J}_\lambda^c=\{(u,v)\in X\times X:\mathcal{J}_\lambda(u,v)\leq c\}$, $\mathscr{U}$ is a neighborhood of $(\varphi_\lambda,\psi_\lambda)\in X\times X$, $\mathsf{G}=\mathbb{Z}$ or $\mathbb{R}$ is the coefficient group of $\mathsf{H}_q$.
 With no loss of generality, we assume $(\varphi_\lambda,\psi_\lambda)$ is the unique critical point of $\mathcal{J}_\lambda$.
 Since $(\varphi_\lambda,\psi_\lambda)$ is a local minimum critical point,  we easily get
 $$\mathsf{C}_q(\mathcal{J}_\lambda,(\varphi_\lambda,\psi_\lambda))=\delta_{q0}\mathsf{G}. $$
  By Lemma \ref{PS},
 $\mathcal{J}_\lambda$ satisfies the Palais-Smale condition. Then applying (\cite{Chang93}, Chapter 2, Theorem 3.2) and Theorem \ref{degree-system}, we obtain
 \bna
 0=\deg{(\mathcal{F},B_R,(0,0))}&=&\deg(D\mathcal{J}_\lambda,B_R,(0,0))\\
 &=&\sum_{q=0}^\infty (-1)^q{\rm rank}\, \mathsf{C}_q\le(\mathcal{J}_\lambda,(\varphi_\lambda,\psi_\lambda)\ri)\\
 &=&1,
 \ena
 provided that $R>\|\varphi_\lambda\|_{L^\infty(V)}+\|\psi_\lambda\|_{L^\infty(V)}$.
 This is impossible, and thus $\mathcal{J}_\lambda$ must have another critical point,
 as we desired.
 $\hfill\Box$\\



\end{document}